 \documentclass[jip]{degruyter-journal-a}      

\usepackage{multirow}

\title{Differential evolution algorithm of solving an inverse problem for the spatial Solow mathematical model}
\headlinetitle{Differential evolution algorithm for spatial Solow model}

\lastnameone{Kabanikhin}
\firstnameone{Sergey}
\nameshortone{S.~Kabanikhin}
\addressone{Institute of Computational Mathematics and Mathematical Geophysics SB RAS, Prosp. Akad. Lavrentyeva~6, Novosibirsk State University, Pirogova str.~2, 630090 Novosibirsk}
\countryone{Russia}
\emailone{kabanikhin@sscc.ru}

\lastnametwo{Krivorotko}
\firstnametwo{Olga}
\nameshorttwo{O.~Krivorotko}
\addresstwo{Institute of Computational Mathematics and Mathematical Geophysics SB RAS, Prosp. Akad. Lavrentyeva~6, Novosibirsk State University, Pirogova str.~2, 630090 Novosibirsk}
\countrytwo{Russia}
\emailtwo{olga.krivorotko@sscc.ru}

\lastnamethree{Bektemessov}
\firstnamethree{Maktagali}
\nameshortthree{M.~Bektemessov}
\addressthree{Abai Kazakh National Pedagogical University, Dostyk ave.~13, 050010 Almaty}
\countrythree{Republic of Kazakhstan}
\emailthree{maktagali@mail.ru}

\lastnamefour{Bektemessov}
\firstnamefour{Zholaman}
\nameshortfour{Zh.~Bektemessov}
\addressfour{Al-Farabi Kazakh National University, Kazybek Bi st.~30, 050040 Almaty}
\countryfour{Republic of Kazakhstan}
\emailfour{jolaman252@gmail.com}

\lastnamefive{Zhang}
\firstnamefive{Shuhua}
\nameshortfive{Sh. Zhang}
\addressfive{Tianjin University of Finance and Economics, Tianjin}
\countryfive{China}
\emailfive{shuhua55@126.com}

\headlineauthor{S. Kabanikhin, O. Krivorotko, M. Bektemessov and others}

\abstract{The differential evolution algorithm is applied to solve the optimization problem to reconstruct the production function (inverse problem) for the spatial Solow mathematical model using additional measurements of the gross domestic product for the fixed points. Since the inverse problem is ill-posed the regularized differential evolution is applied. For getting the optimized solution of the inverse problem the differential evolution algorithm is paralleled to 32 kernels. Numerical results for different technological levels and errors in measured data are presented and discussed.}



\researchsupported{This work is supported by the Russian Science Foundation (grant No. 18-71-10044), i.e. numerical investigation of spatial Solow mathematical model, and by the grant of the Ministry of Education and Science of the Republic of Kazakhstan (project No.~AP05134121 ''Numerical methods of identifiability of inverse and ill-posed problems of natural sciences''), i.e. inverse problem statement (Section 2).}

\acknowledgments{Authors thank to Professor Daniyar Nurseitov for problem statement and fruitfull discussions for concerning Solow mathematical model and thank to Dr.~Igor Chernykh for help in parallel realization of DE algorithm.}

\begin{document}

\newcommand*\Laplace{\mathop{}\!\mathbin\bigtriangleup}

\section{Introduction}\label{introduction}
As the Solow growth model~\cite{Solow_1956} was build using production function (naturally Cobb-Douglas production function), law of motion for the stock of capital and saving/investment function, model can be easily extended to include a households problem (the Ramsey-Cass-Koopmans model). Usually the interest of the Solow model is that it perpetual growth, that can be obtained using balanced growth path and technological progress over time. Output per worker can grow only as long as capital per worker grows and the key to constant growth is the existence of non-diminishing marginal product of capital. Another way of perpetual economical growth is letting technological progress change in model, it means allowing technological parameter to grow exogenously over time.

As it has been written in~\cite{Enkhbat_2010}, the Solow growth model was considered with assumptions as concave homogeneous production function (instead of Cobb-Douglas production function), exponentially growing labor and constant saving function. In addition, they considered the per capita consumption maximization problem subject to economic equilibria. Authors considered two cases when production function is logistic and the labor grow exponentially and when both of them are logistic, to reduce them to one variable parametric maximization problem. After making sure that production function is nonconvex and satisfies the Lipschitz conditions, authors solved nonconvex optimization problem by global optimization techniques, considered on a sufficiently large interval. The solution can be found by the method of piecewise linear function~\cite{Horst_1995}.

In the article~\cite{Smirnov_2018} by Smirnov and Wang, the work of Ryuzo Sato~\cite{Sato_1980}, devoted to the development of economic growth models within the framework of the Lie group theory, was extended to a new growth model based on the assumption of logistic growth by using the Solow economic growth model as a starting point. Authors claimed that the Cobb-Douglas function can no longer adequately describe the growth of the economy over a long-run, it was aimed to develop a new mathematical paradigm that can be used to study the current state of economy and to replace neoclassical growth model in the sense of Sato representing exponential growth with a logistic growth. Also they used the new ``S-shaped'' production function, the consequence of logistic growth in factors, to solve maximization problem of profit under condition of perfect competition, using the same arguments of subject to relevant changes by assuming that the revenue of the firm from sales is determined.

The logistic growth in other words can be described by spatial Solow model and that was used in~\cite{Engbers_2014}, where they did identification of production function using (noisy) data that is an ill-posed inverse problem, using non-parametric approach and applied Tikhonov regularization to stabilize the computations. As there is no clear choice which production function will fit the situation best, it was proposed to identify production function from data about the capital distribution of some spatial economy, further they obtained the following optimization problem that has to be minimized. To solve the minimization problem authors applied the gradient descent algorithm. As the objective function was Freshet differentiable, they used the directional derivatives of the Langrangian and to find the minimum of the functional the simple steepest descent method and a backtracking line search method~\cite{Nocedal_1999} were used. So they reconstructed the production function to a spatial Solow model with the different noise levels and different technology terms, when it is constant and space-dependent.

We use the spatial Solow mathematical model as in~\cite{Engbers_2014} and investigate the inverse problem for its using stochastic approach for global optimization.

The paper is organized as follows. Section~2 presents the derivation of the original Solow mathematical model described by ordinary differential equation and the statement of inverse problem for the Solow model described by partial differential equation. The formulation of an inverse problem as the optimization problem and numerical algorithm for solving inverse problem is presented in Section~3. The results of numerical calculations for spatial Solow model are presented and discussed in Section~4. Conclusions are given in Section~5, followed by a list of references.

\section{Statement of the problem}
In this Section the derivation of neoclassical Solow mathematical model for ordinary differential equation is demonstrated at subsection~\ref{ODE_Solow}. Based on that derivation the statement of the spatial Solow model is considered in subsection~\ref{spatial_Solow} and the inverse problem statement for spatial Solow model is formulated at subsection~\ref{inverse_Solow}. 

\subsection{Solow mathematical model for ODE}\label{ODE_Solow}
Neoclassical economical Solow model describes evolution of gross output -- $Y(t)$, using next (due to such) indicators as: used labor resources -- $L(t)$, saving capital -- $K(t)$ and technological progress -- $A(t)$. And since the output parameter of the model should be a stable indicator of a productive economy, then the gross domestic product (GDP) is taken, which is a macroeconomic index reflecting the market value of all final goods and services produced during the year in the state~\cite{McConnellBrueFlynn}. A mathematical notation connecting these variables is $Y(t)=A(t)Q(K(t),L(t))$,
where $Q$ represents production function. It is assumed that the production function is homogeneous, which means $Q(\alpha K(t),\alpha L(t))=\alpha Q(K(t),L(t))$.
Also it can be noted that the production function satisfies the following condition
\begin{eqnarray*}
	Q( 0, L(t))=0= Q(K(t),0).
\end{eqnarray*}
It is worth saying that the description of the development of any economy only due to the absolute value of any gross output is useless, it is hard to say whether the economy is doing well or not. Simon Kuznets, one of the architects of the US national accounting system, the man who first introduced the concept of GDP in 1934, warned against identifying GDP growth with increasing economic or social welfare.What we are interested in is the rate of economic growth \cite{Kuznets_1934, Kuznets_1941}. Therefore, we consider the rate of change in capital, which looks like 
\begin{eqnarray}\label{eq:change_capital}
	\dfrac{dK}{dt}=Y(t)-C(t)-\delta K(t).
\end{eqnarray}
It means the change in fixed capital stock negatively depends on the volume of consumption $C(t)$  and on the amount of depreciation that is supposed to occur with the rate $\delta$. Moreover, we assume that the difference in production and consumption persists for each period of time, namely
\begin{eqnarray}\label{eq:prod_cons}
	Y(t)-C(t)=sY(t).
\end{eqnarray}
Then using that and inserting~(\ref{eq:prod_cons}) into~(\ref{eq:change_capital}), we have the following
\begin{eqnarray}\label{eq:gross_output_sub}
	\dfrac{dK}{dt}=sA(t)Q(K(t),L(t))-\delta K(t).
\end{eqnarray}
Next, we introduce a new variable, namely $k(t)=\dfrac{K(t)}{L(t)}$, the capital per capita. Then it turns out, using the homogeneous of function $Q$ we can write the following 
\begin{eqnarray*}
	q(k(t))=\frac{1}{L(t)}Q(K(t),L(t))=Q\left(\frac{K(t)}{L(t)},1\right) 
\end{eqnarray*}
and calculate 
\begin{eqnarray*}
	\frac{dk(t)}{dt}=\frac{d}{dt}\left(\frac{K(t)}{L(t)}\right)=\frac{\frac{dK(t)}{dt}}{L(t)}-n\frac{K(t)}{L(t)},
\end{eqnarray*}
where $n=\frac{\frac{dL(t)}{dt}}{L(t)}$ denotes a constant growth rate of labor costs (labor intensity). With these designations and abbreviations, we can rewrite~(\ref{eq:gross_output_sub}) as follows 
\begin{eqnarray}\label{eq:gross_output_sub_rewrite}
	\dfrac{dk(t)}{dt}=sA(t)q(k(t))-(\delta+n)k,
\end{eqnarray}
which is the basic equation for spatial structured Solow model~\cite{Solow_1956, Krugman_1991, Mossay_2003}.

It is worthy to clarify that we are interested in the change in capital for a work unit (that is, an employee) - the capital-labor ratio, or more precisely, the situation where the capital per work unit reaches its steady state. To do this, consider a stationary solution of equation ~(\ref{eq:gross_output_sub_rewrite})
\begin{eqnarray*}
	0=sA(t)q(k(t))-(\delta+n)k.
\end{eqnarray*}
If we assume that $A(t)=1$, then it means that there is no technological progress at all. Then there are only three variables describing the capital-labor ratio: saving rate - $s$, depreciation rate - $\delta$ and the rate of population growth or unit of labor used - $n$. Consequently, capital intensity will increase (grow) if 
\begin{eqnarray*}
	sq(k(t))>(\delta+n)k
\end{eqnarray*}
and decrease (fall) otherwise. Thus, if the capital ratio is a constant number, then economy tends to its steady state $k_E$, i.e. there are enough savings to cover the costs associated with population growth and the amount of capital lost due to depreciation. Moreover, the economic growth rate in steady state equals to rate of population growth (i.e. $n$). Further, we assume that parameters such as population growth rate and depreciation coefficient are always constant, then the only variable affecting the model is the savings rate - $s$. It is also assumed that when saving changes from $s_0$ to $s_1$ at ($s_1>s_0$), the function shows a sharp rise, and then the steady state increases from $k_0$ to $k_1$. It is good for economy for a short period of time, because economic growth occurs faster, but in the long run the economy will tend to a new steady state and then the economic growth rate will again be equal to $n$. So $n$ is not only constant, but also equals zero, since the population does not change at all. The rate of savings over a large time interval, in turn, does not have any effect either on the rate of economic growth. The only option to obtain economic growth is a technological progress~\cite{Engbers_2009}. Thus, if the parameter as $n$ does not have any effect on model - $n=0$, noting that we set $s=1$ for simplification, then the equation~(\ref{eq:gross_output_sub_rewrite}) should be rewritten as
\begin{eqnarray}\label{eq:basic_rewr}
	\dfrac{dk(t)}{dt}=A(t)q(k(t))-\delta k.
\end{eqnarray}

\subsection{The spatial Solow model}\label{spatial_Solow}
Consider the scaled initial-boundary value problem for the mathematical model described dynamic of the capital stock held by the representative household located at $x$ at date $t$~\cite{Camacho_2008, Engbers_2014}. Then the mathematical model~(\ref{eq:basic_rewr}) with adding initial and boundary conditions is rewritten as follows:
\begin{eqnarray}\label{eq:spartial_Solow_m}
\left\{\begin{array}{ll}
	\dfrac{\partial k}{\partial t} - d \Laplace k(x,t) = g(k, x, t), \quad & x\in \varOmega , t\in [0,T],\\
	k(x,0) = k_0(x) > 0, & x\in \varOmega,\\
	\nabla k\cdot n = 0, & \text{on}\,\,\partial \varOmega\times [0,T].
\end{array}\right.
\end{eqnarray}
Here $d=\frac{1}{\delta L^2}$ is a scaled coefficient, $\delta$ is the depreciation rate, $g(k,x,t) = \frac{A(x,t)}{\delta}q(k) - k$, $A(x,t)$ denotes the technological level at $x$ and time $t$. The standard neoclassical production function is assumed to be non-negative, increasing and concave, and verifies the Inada conditions, that is,
\begin{eqnarray*}
	\lim_{k\to 0}q^\prime (k) = +\infty, \quad \lim_{k\to \infty}q^\prime (k) = 0,\quad q(0)=0.
\end{eqnarray*}
We will depart from the assumptions with respect to concavity in particular around zero as well as the first Inada condition and allow for general convex-concave production functions, an example being~\cite{Engbers_2014}
\begin{eqnarray}\label{eq:prouction_func}
	q(k)=\dfrac{\alpha_1 k^p}{1+\alpha_2 k^p},\quad \alpha_1,\alpha_2 \geq 0, p > 1.
\end{eqnarray}
Such examples of $q$ are of particular interest, because they are related to the potential existence of poverty traps. Define the set of admissible production functions
\begin{eqnarray*}
\begin{array}{cc}
Q_{\mbox{adm}} = \left\{ q\in H^1(0,K)\ |\ q(0)=0, 0\leq q^\prime(k)\leq q^\prime_{\max}\,\, \text{for}\,\, k\in(0,K),\right.\\
\qquad\qquad\qquad\quad \left. q^\prime(k) = 0\,\, \text{else} \right\},
\end{array}
\end{eqnarray*}
where $q^\prime_{\max}$ being a fixed constant, which can be understood as the maximal growth that an economy is capable of.

The technological level $A(x,t)$ is determined via a diffusion equation of the form
\begin{eqnarray}\label{eq:tech_level_A}
\left\{\begin{array}{ll}
	\dfrac{\partial A}{\partial t} - \Laplace A = Ag_A, \quad & x\in \varOmega , t\in [0,T],\\
	A(x,0) = A_0(x), & x\in \varOmega,\\
	\frac{\partial A}{\partial x} = 0, & \text{on}\,\,\partial \varOmega\times [0,T],
\end{array}\right.
\end{eqnarray}
with $g_A$ being either constant, a function depending only on space or a function depending on space as well as on time.

The Neumann boundary condition in problem~(\ref{eq:spartial_Solow_m}) represents no capital flow through the boundary and thereby a closed economy.

In paper~\cite{Engbers_2014} authors proved a well-posedness of direct problem~(\ref{eq:spartial_Solow_m}) at space $L^2([0,T], H^1(\varOmega))\cap H^1([0,T], H^{-1}(\varOmega))$ if $k_0\in L^\infty(\varOmega)$, $q\in Q_{\mbox{adm}}$ and $A\in C(\varOmega\times [0,T])$. A more detailed analysis of this model can be found in~\cite{Capasso_2010}.

\subsection{Inverse problem statement for spatial Solow model}\label{inverse_Solow}
The choice of the production function is crucial for an economic model, as its shape will greatly influence the capital distribution. In general, data about the economic situation, such as the gross domestic product (GDP), of different regions and different countries are readily available. Suppose, that we have additional information about GDP of some spatial economy at fixed space and time points:
\begin{eqnarray}\label{eq:ip_data}
	k(x_m, t_j) = f_{mj} + \varepsilon_{mj},\; x_m\in \varOmega, t_j\in [0,T], \, m=1,\ldots,M, j=1,\ldots, N.
\end{eqnarray}
Here $\varepsilon_{mj}$ are Gaussian noise in measurements.

The \textit{inverse problem}~(\ref{eq:spartial_Solow_m}),~(\ref{eq:ip_data}) consists in identification of production function~(\ref{eq:prouction_func}) (or identification of parameters $\alpha_1,\alpha_2,p$) of initial-boundary value problem~(\ref{eq:spartial_Solow_m}) using additional measurements~(\ref{eq:ip_data}). It means that we have the nonlinear parameter-to-solution map $A:\ q\in Q_{\text{adm}} \mapsto f^\varepsilon\in E^{MN}$ mapping the production function $q$ to the respective capital
distribution $f^\varepsilon = \{f_{mj} + \varepsilon_{mj}\}_{\substack{m=1,\ldots,M,\\ j=1,\ldots, N }}$, i.e. $A(q) = f^\varepsilon$.
Here $E$ is an Euclidean space of measurements.

The inverse problem~(\ref{eq:spartial_Solow_m}),~(\ref{eq:ip_data}) is ill-posed~\cite{KashtanovaVNKabanikhin2011}, i.e. the solution $q(k)$ is non-unique and can be unstable~\cite{Engbers_2014}. That is we apply the regularization technique described in Section~\ref{sec:optimization}.

\section{Optimization problem and numerical algorithm}\label{sec:optimization}
Reduce our inverse problem~(\ref{eq:spartial_Solow_m}),~(\ref{eq:ip_data}) to an optimization problem that consists in minimization of the misfit function
\begin{eqnarray}\label{misfit_func}
	J(q) = \Vert A(q) - f^\varepsilon\Vert^2_{L^2_\chi(\varOmega\times [0,T])} := \int\limits_0^T\int\limits_\varOmega \chi(x,t) (A(q) - f^\varepsilon)^2\, dxdt.
\end{eqnarray}
Here $\chi (x,t)$ is a characteristic function of incomplete measurements~(\ref{eq:ip_data}). In our case the misfit function~(\ref{misfit_func}) has the form:
\[J(q) = \dfrac{1}{N M}\sum\limits_{j=1}^N\sum\limits_{m=1}^M (k(x_m,t_j;q) - f_{mj}^\varepsilon)^2.\]

Optimization problem can be solved by various methods such as gradient approaches, stochastic methods, etc~\cite{KOI_KSI_arXiv2019}. The misfit function~(\ref{misfit_func}) has a lot of local minimums due to ill-posedness of inverse problem~(\ref{eq:spartial_Solow_m}),~(\ref{eq:ip_data}). In paper~\cite{Engbers_2014} authors applied the Tikhonov regularization approach based on gradient method with Tikhonov regularization term. The main weaknesses of this approach are the difficulty of choosing the regularization parameter and the dependence of the convergence of the gradient method on the choice of the initial approximation (local convergence). We choose the stochastic algorithm of global optimization based on solving more simple evolutionary problems from biology named differential evolution algorithm~\cite{Storn_Price_1995}.

\subsection{Differential evolution algorithm}\label{DE_alg}
Differential evolution algorithm (DE), a class of evolutionary algorithms, was introduced
by Storn and Price at 1995~\cite{Storn_Price_1995, Storn_Price_1997, Price_2005, Storn_2008} for solving a polynomial fitting problem. The algorithm is generally called as a very simple but very powerful population-based meta-heuristic algorithm~\cite{Qing_2009}.
The algorithm is generally characterized by the features of simplicity, effectiveness and robustness. Also, it is easy-to-use, and it requires few controlling parameters, and it
has fast convergence characteristic~\cite{Storn_Price_1997}. Due to these advantages, it presents a wide range of implementation examples in different areas such as acoustics, biology, material science, mechanic, medical imaging, optic, mathematics, physics, seismology, economics etc. More details and examples about the implementation of DE to solve various problems are given in~\cite{Qing_2009}.
Even though previous comprehensive studies over real-world problems have shown that DE performs better in terms of convergence rate and robustness~\cite{Sambarta2009_DE_convergence, Hahn_1963} than the other evolutionary algorithms such as genetic algorithm, particle swarm optimization~\cite{PSO_1995}, simulated annealing~\cite{Kirkpatrick_1983}, etc.

An algorithm of differential evolution is follows:
\begin{enumerate}
\item[1.] \textit{Initialization}. Create an initial population of target vectors of parameters $q_{i,G} = \left( q_{i,G}^1, q_{i,G}^2, q_{i,G}^3 \right)$, $i=1,\ldots,Np$, where $Np$ is the population size, $G$ denotes current generation. Here $q_{i,G}^1 = \alpha_{1_{i,G}}$, $q_{i,G}^2 = \alpha_{2_{i,G}}$, $q_{i,G}^3 = p_{i,G}$. The algorithm is initialized by a randomly created population within a predefined search space considering the upper (index $u$) and lower (index $l$) bounds of each parameter $q_{i,G}^j \in [q_{l}^j, q_{u}^j]$, $j=1,2,3$.
\item[2.] \textit{Choose stopping criteria}. Set the stopping parameter $\varepsilon_{stop}$ for misfit function  and maximum number of iterations $G_{\max}$. If $J(q_{i,G}) < \varepsilon_{stop}$ for some $i=1,\ldots, Np$ or $G=G_{\max}$ then stop iterations and choose $i$ with minimum value of misfit function $J(q_{i,G})$. Otherwise go to step 3.
\item[3.] \textit{Mutation}. At each iteration, the algorithm generates a new generation of vectors, randomly combining vectors form the previous generation. For each new generation ($G+1$) of a vector from a given target vector $q_{i}$ from the old generation ($G$) algorithm randomly selects three vectors $q_{r_1,G}$, $q_{r_2,G}$ and $q_{r_3,G}$ such that $i, r_1, r_2, r_3$ are distinct and creates a donor vector 
\[v_{i,G+1} = q_{r_1,G} + F(q_{r_2,G} - q_{r_3,G}), \quad F\in [0,2] \; \text{is a differential weight.}\]
\item[4.] \textit{Crossover (recombination)}. Create the trial vector $u_{i,G}$ from the elements of the target vector $q_{i,G}$ and donor vector $v_{i,G+1}$ with probability $Cr\in [0,1]$ using formula:
\begin{eqnarray*}
u_{i,G+1}^j = \left\{\begin{array}{ll}
v_{i,G+1}^j, \quad & \text{if}\;\; \mbox{rand}_{i,j}\leq Cr\;\; \text{or}\;\; j=j_{\mbox{rand}},\\
q_{i,G}^j, & \text{otherwise}
\end{array}\right., \; j=1,2,3.
\end{eqnarray*}
Here $\mbox{rand}_{i,j}$ represents a uniformly distributed random variable in the range of $[0,1)$, $j_{\mbox{rand}}$ is a randomly chosen integer in the range $[1, 3]$ to
provide that the trial vector does not duplicate the target vector.
\item[5.] \textit{Selection}. The  vector obtained after crossover is the test vector. If it is better than the base vector, then in the new generation the base vector is replaced by trial one, otherwise the base vector is stored in the new generation. Choose the next generation as follows:
\begin{eqnarray*}
q_{i,G+1} = \left\{\begin{array}{ll}
u_{i,G+1}, \quad & J(u_{i,G+1})\leq J(q_{i,G}),\\
q_{i,G}, & \text{otherwise}.
\end{array}\right.,\quad i=1,\ldots, Np
\end{eqnarray*}
and go to step 2 till $G+1<G_{max}$.
\end{enumerate}

\section{Numerical experiments}
%
%
%
%
%
%
%

We will show some numerical results of the inverse problem~(\ref{eq:spartial_Solow_m}),~(\ref{eq:ip_data}) using DE algorithm described is Section~\ref{DE_alg}. We start by giving some details about the simulated dataset used for the calculations and then show the identification results for a constant technological level $A$~(Section~\ref{constantA}) and for a space-dependent technology term~$A(x)$ (Section~\ref{spacedepA}).

\subsection{Simulated dataset}\label{simulated_dataset}
Consider the modelling scaled domain $[0,L]\times [0,T]$ with $L=50$ and $T=150$ (here $L$ can numerated regions with different GDP and $T$ described time in years). After nondimensionalization we get new computational domain $[0,1]\times [0,\delta T]$ where mathematical problem~(\ref{eq:spartial_Solow_m}) is formulated. We put $\delta = 0.05$. We set an equidistant grid with $N_x=26$ nodal points in space and $N_t=251$ nodal points in time, which leads to a spatial-step size $h_x = 1/N_x = 0.04$ and a time-step size $h_t =\delta T/N_t = 0.03$. The classic second-order difference approximation has been used to discretize the diffusion. The time derivative is approximated by backward difference of the first-order.

We put an initial condition $k_0(x)$ as a piece-wise function on the interval $[0,1]$:
\begin{eqnarray*}
k_0(x) = \left\{\begin{array}{ll}
0, & x\in [0,0.3),\\
25(x-0.3), & x\in[0.3, 0.7],\\
10, &x\in(0.7,1].
\end{array}\right.
\end{eqnarray*}
We obtain the synthetic data $f_{mj}$ from~(\ref{eq:ip_data}) for different $M$ and $N$ by solving the direct problem~(\ref{eq:spartial_Solow_m}) with the production function
\[q_{ex}(k) = \dfrac{0.0005 k^4}{1+0.0005 k^4}\]
presented at figure~\ref{ris:func_q} (left) and two types of technological terms $A(x,t)$ (see below). Measurements are uniformly distributed in space on $[0,1]$ and time on $[\delta T/2, \delta T]$ (see example for $M=5$, $N=6$ at figure~\ref{ris:func_q} right).
\begin{figure}[!ht]
\begin{center}
\includegraphics[width=1\linewidth]{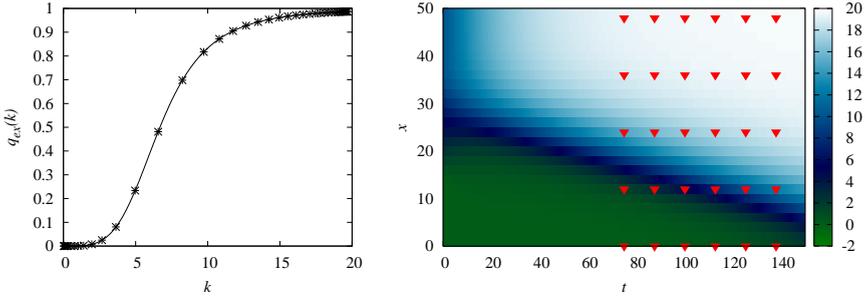}
\caption{The exact production function $q_{ex}(k)$ (left) and map of direct problem solution $k(x,t;q_{ex})$ with points $m=1,\ldots,M$, $j=1,\ldots,N$ of measurements~(\ref{eq:ip_data}) for $M=5$, $N=6$ (right).}
\label{ris:func_q}
\end{center}
\vspace*{-2mm}
\end{figure}
Then we add the Gaussian noise to inverse problem data~(\ref{eq:ip_data}) as follows
\[f_{mj}^\varepsilon = f_{mj} + \varepsilon f_{mj} \xi_{mj},\quad m=1,\ldots,M, j=1,\ldots,N.\]
Here $\xi_{mj} \sim N_{0,1}$ is a normally distributed modeled random variables with zero mean and unit dispersion, $\varepsilon$ is an error level.

For DE algorithm we put population size $Np=100$ and we choose parameters $F=0.7$ and $Cr=0.9$ as the best combination for convergence features for the algorithm~\cite{BALKAYA2013160}. We set maximum number of iterations $G_{\max} = 5000$ and $\varepsilon_{stop} = 10^{-4}$. For getting the optimized solution of the inverse problem we launch the DE algorithm 1000 times for all decribed numerical calculations using the cluster NKS-30T in the Siberian Supercomputer Center in the Institute of Computational Mathematics and Mathematical Geophysics of the SB RAS and then take the arithmetic average. 

\subsection{Numerical results for constant technological term $A(x,t)=1$}\label{constantA}
\begin{figure}[!ht]
\begin{center}
\includegraphics[width=1\linewidth]{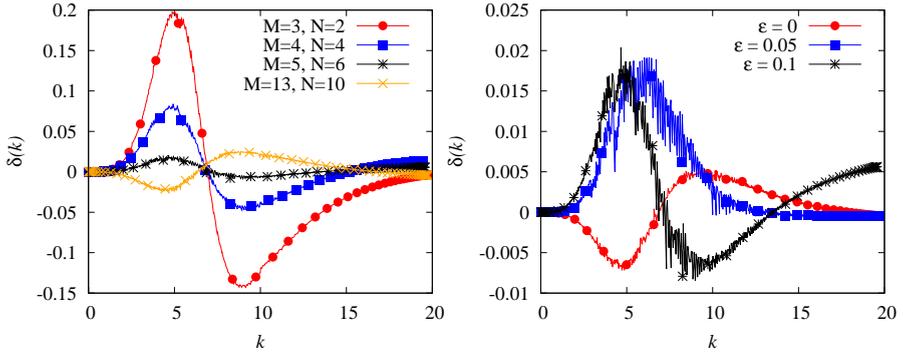}
\caption{The difference $\delta (k)$ of exact and approximate solutions for different points of measurements $M$ and $N$ for fixed error level in data~(\ref{eq:ip_data}) $\varepsilon=0.1$ (left). The difference $\delta (k)$ of exact and approximate solutions for different noise levels $\varepsilon=0, 0.05, 0.1$ in measurements for $M=5$, $N=6$ (right).}
\label{ris:func_q_difMN}
\end{center}
\vspace*{-2mm}
\end{figure}

We solve optimization problem~(\ref{misfit_func}) with constant technological term $A(x,t) = 1$ using DE described in Section~\ref{DE_alg}. For $\varepsilon =0.1$ in data~(\ref{eq:ip_data}) we get the inverse problem solution $q_\varepsilon(k)$ for four variants of $M$ and $N$. Figure~\ref{ris:func_q_difMN} (left) demonstrates the difference $\delta(k) = q_{ex}(k)-q_\varepsilon(k)$ of exact and calculated solutions of inverse problem with four variants of measurements. Table~\ref{tab:tab1} shows that the $M=5$, $N=6$ is sufficient for reconstruction of production function with necessary accuracy in relative error $\rho_\varepsilon = \Vert k(\cdot,\cdot;q_{ex}) - k(\cdot,\cdot;q_\varepsilon)) \Vert_{L^2}/\Vert k(\cdot,\cdot;q_{ex}) \Vert_{L^2}$. The smaller number of measurement points the greater the difference $\delta(k)$ (see figure~\ref{ris:func_q_difMN} left).

\begin{table}[ht!]
\caption{Relative errors and value of the misfit function $J(q_\varepsilon)$ for different number of measurements~(\ref{eq:ip_data}) with error level $\varepsilon = 0.1$ and constant technological level $A(x,t)=1$.} \label{tab:tab1}
{\footnotesize
\begin{tabular}{cccc} \hline
Values of $M$ and $N$ & $\max|\delta(k)|$ & $\rho_\varepsilon$ & $J(q_\varepsilon)$ \\
\hline
$M=3$, $N=2$ & 0.171 & 0.024 & 0.102 \\
$M=4$, $N=4$ & 0.053 & 0.015 & 0.479 \\
$M=5$, $N=6$ & 0.02 & 0.004 & 0.293\\
$M=13$, $N=10$ & 0.01 & 0.006 & 0.199 \\
\hline
\end{tabular}}
\vspace*{-2.5mm}
\end{table}

For $\varepsilon =0$, $0.05$ and $0.1$ in data~(\ref{eq:ip_data}), $M=5$, $N=6$, we get the inverse problem solution $q_\varepsilon(k)$ (the differences $\delta(k)$ are plotted on figure~\ref{ris:func_q_difMN} right).Table~\ref{tab:tab2} shows the reconstructed parameters $\alpha_1$, $\alpha_2$ and $p$ in function~(\ref{eq:prouction_func}) for different error level in measured data. If we have noise free data of inverse problem then the difference $\delta(k)$ is close to zero. It means that reconstruction of parameters $\alpha_1, \alpha_2, p$ is close to the tested ones (the maximum of the absolute difference $\delta(k)$ is equal to 0.005 as given in table~\ref{tab:tab2}). Note, that maximum absolute error of inverse problem solutions for $M=5$, $N=6$ is less than 2\%, i.e. $\max|\delta(k)| \leq 0.02$ for maximum error level in inverse problem data~(\ref{eq:ip_data}) $\varepsilon = 0.1$.
\begin{table}[ht!]
\caption{Reconstructed parameters in function $q_\varepsilon(k)$ for different error levels $\varepsilon = 0, 0.05, 0.1$ in measurements~(\ref{eq:ip_data}) for $M=5$, $N=6$ and constant technological term $A(x,t) = 1$.} \label{tab:tab2}
{\footnotesize
\begin{tabular}{ccccc} \hline
\multirow{2}{*}{Parameters} & \multirow{2}{*}{Exact values} & \multicolumn{3}{c}{Error level in inverse problem data} \\
& & $\varepsilon = 0$ & $\varepsilon = 0.05$ & $\varepsilon = 0.1$ \\
\hline
$\alpha_1$ & 0.0005 & 0.00057  & 0.00048 & 0.00036 \\
$\alpha_2$ & 0.0005 & 0.00057 & 0.00048 & 0.00036 \\
$p$ & 4 & 3.9226 & 4.0202 & 4.1806 \\
\hline
$\max|\delta(k)|$ & & 0.005 & 0.019 & 0.02 \\
\hline
\end{tabular}}
\vspace*{-2.5mm}
\end{table}

\begin{figure}[t]
\begin{center}
\includegraphics[width=1\linewidth]{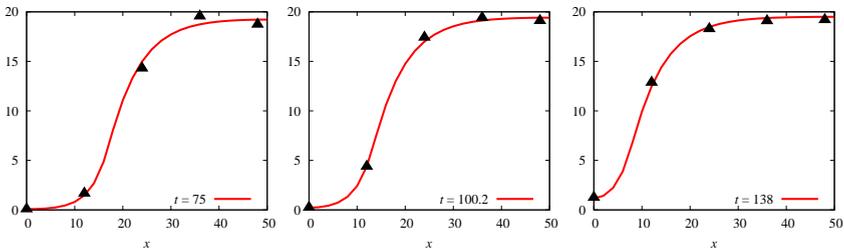}
\caption{The solution $k(x,t_j;q_\varepsilon)$ of the direct problem~(\ref{eq:spartial_Solow_m}) for reconstructed $q_\varepsilon(k)$ with error level in measurements $\varepsilon = 0.1$ and points of measurement are $M=5$, $N=6$ for constant technological term. Here $t_1=75$, $t_3=100.2$, $t_6=138$.}
\label{ris:direct_pr_solution}
\end{center}
\vspace*{-2mm}
\end{figure}

The solution $k(x,t;q_\varepsilon)$ of spatial Solow mathematical model for reconstructed $q_\varepsilon(k)$ and measured data~(\ref{eq:ip_data}) for $\varepsilon = 0.1$ is demonstrated on figure~\ref{ris:direct_pr_solution}.
Figure~\ref{ris:direct_pr_solution} shows the compliance of model solution $k(x,t;q_\varepsilon)$ (red line for fixed time point) with measured synthetic noisy data with noise level $\varepsilon = 0.1$ (black triangles).

\subsubsection{Sensitivity analysis of spatial Solow mathematical model.} Investigate the influence of parameters $\alpha_1$, $\alpha_2$ and $p$ to the mathematical model~(\ref{eq:spartial_Solow_m}) namely to the right-hand side
\[g(k, x,t) = \gamma \frac{\alpha_1 k^p}{1+\alpha_2 k^p} - k,\quad \gamma = \frac{A(x,t)}{\delta}.\]
For this function $g$, consider its gradient by parameters:
\[\frac{\partial g}{\partial\alpha_1} = \frac{\gamma k^p}{1+\alpha_2 k^p}, \; \frac{\partial g}{\partial\alpha_2} = -\frac{\gamma\alpha_1 k^{2p}}{(1+\alpha_2 k^p)^2}, \; \frac{\partial g}{\partial p} = -\frac{\gamma\alpha_1\mbox{ln}(k) k^{p}}{(1+\alpha_2 k^p)^2}.\]
\begin{figure}[h!]
\begin{minipage}[h]{0.6\linewidth}
\center{\includegraphics[width=\linewidth]{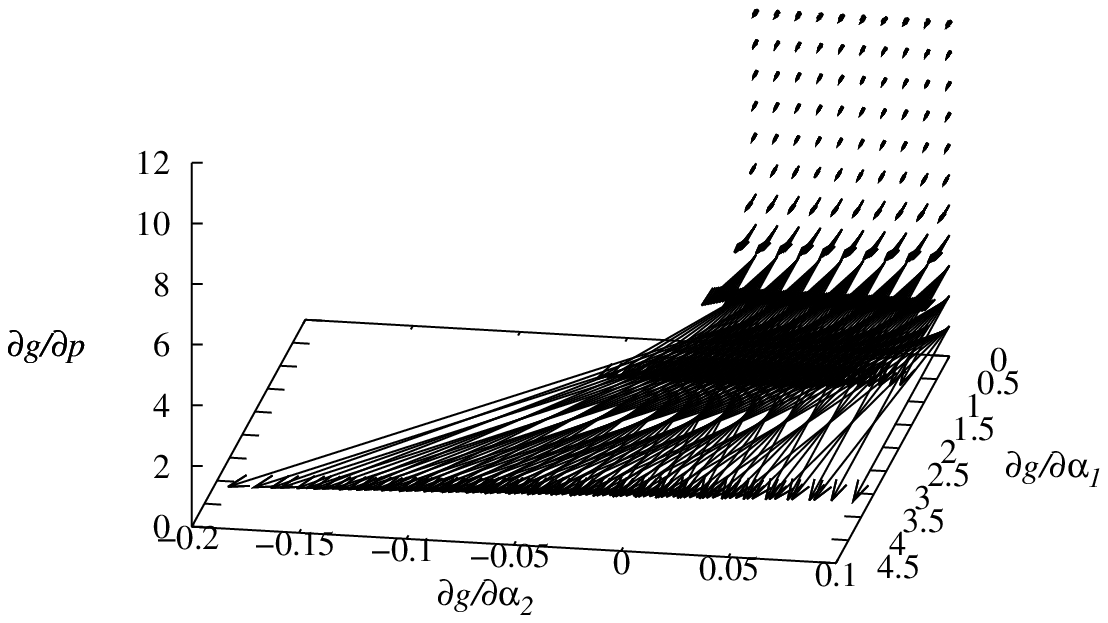} a)}
\end{minipage}
\vspace*{-10mm}
\hfill
\begin{minipage}[h!]{0.6\linewidth}
\center{\includegraphics[width=\linewidth]{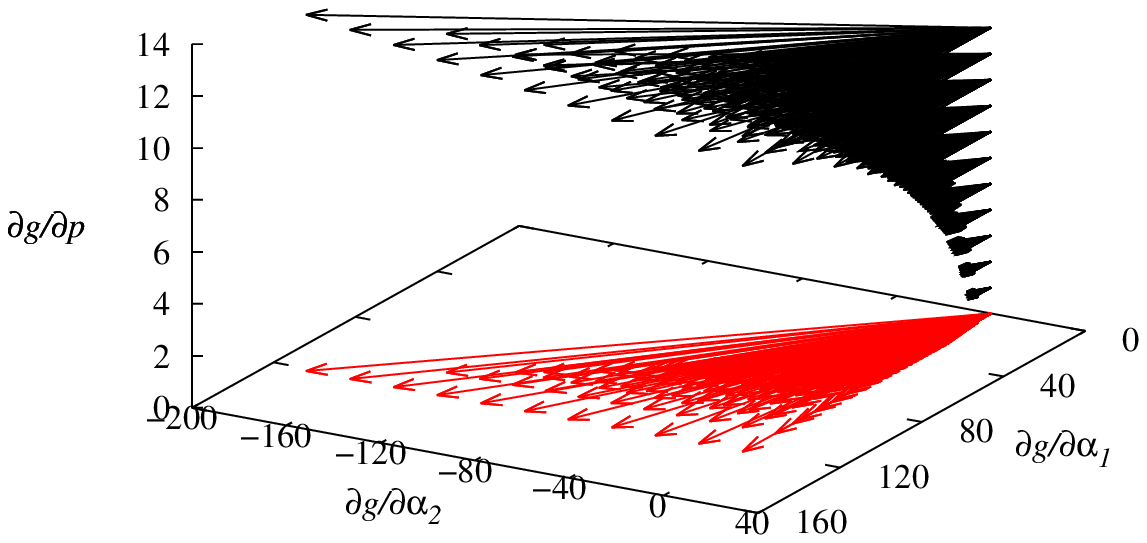} b)}
\end{minipage}
\vspace*{-10mm}
\hfill
\begin{minipage}[h!]{0.6\linewidth}
\center{\includegraphics[width=\linewidth]{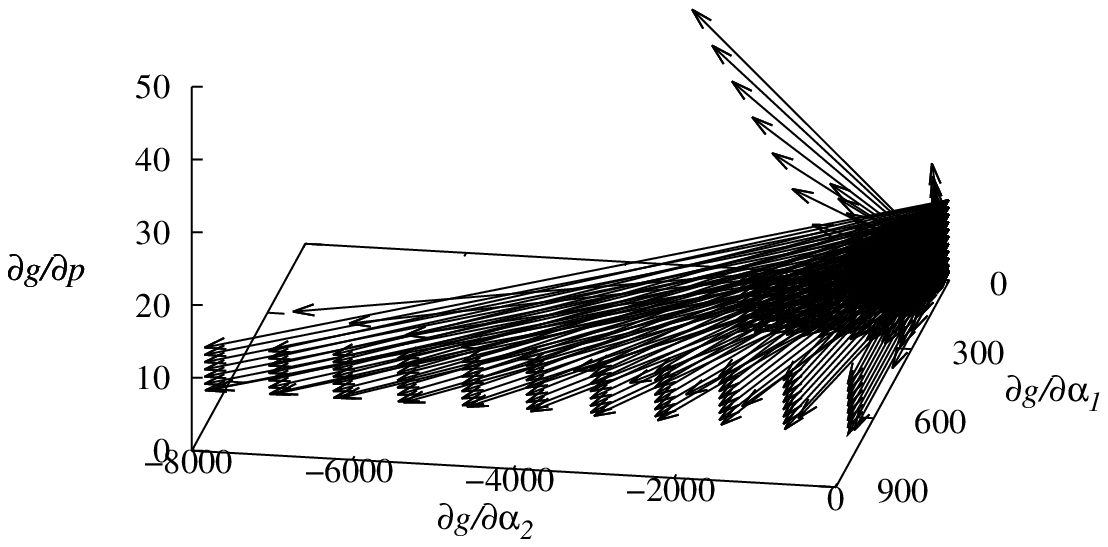} c)}
\end{minipage}
\caption{The gradient field $(\frac{\partial g}{\partial\alpha_1}, \frac{\partial g}{\partial\alpha_2}, \frac{\partial g}{\partial p})$ for fixed values of $k(x,t)$: a) $k(x,t) = 0.5$, b) $k(x,t) = 1.3$ and c) $k(x,t) = 11.25$. Red color arrows show the projection of the gradient field on $(\alpha_1,\alpha_2)$ plane.}
\label{ris:grad_sensitivity}
\end{figure}
For different values of function $k(x,t)$ we construct the gradient field of function $g$. Figure~\ref{ris:grad_sensitivity} shows that the maximum rate of gradient variability for small values of capital stock $k(x,t)$ corresponds to parameters $\alpha_1$ and $\alpha_2$. For bigger values of $k(x,t)$ (figure~\ref{ris:grad_sensitivity} right) and for small values of parameters $\alpha_1$ and $\alpha_2$ gradient grows to the direction of parameter $p$, but when the values of parameters $\alpha_1$ and $\alpha_2$ became bigger the gradient growth turns to parameters $\alpha_1$ and $\alpha_2$ again.

\subsection{Numerical results for space dependent technological level $A(x)$}\label{spacedepA}
\begin{figure}[h!]
\begin{minipage}[h]{0.47\linewidth}
\center{\includegraphics[width=\linewidth]{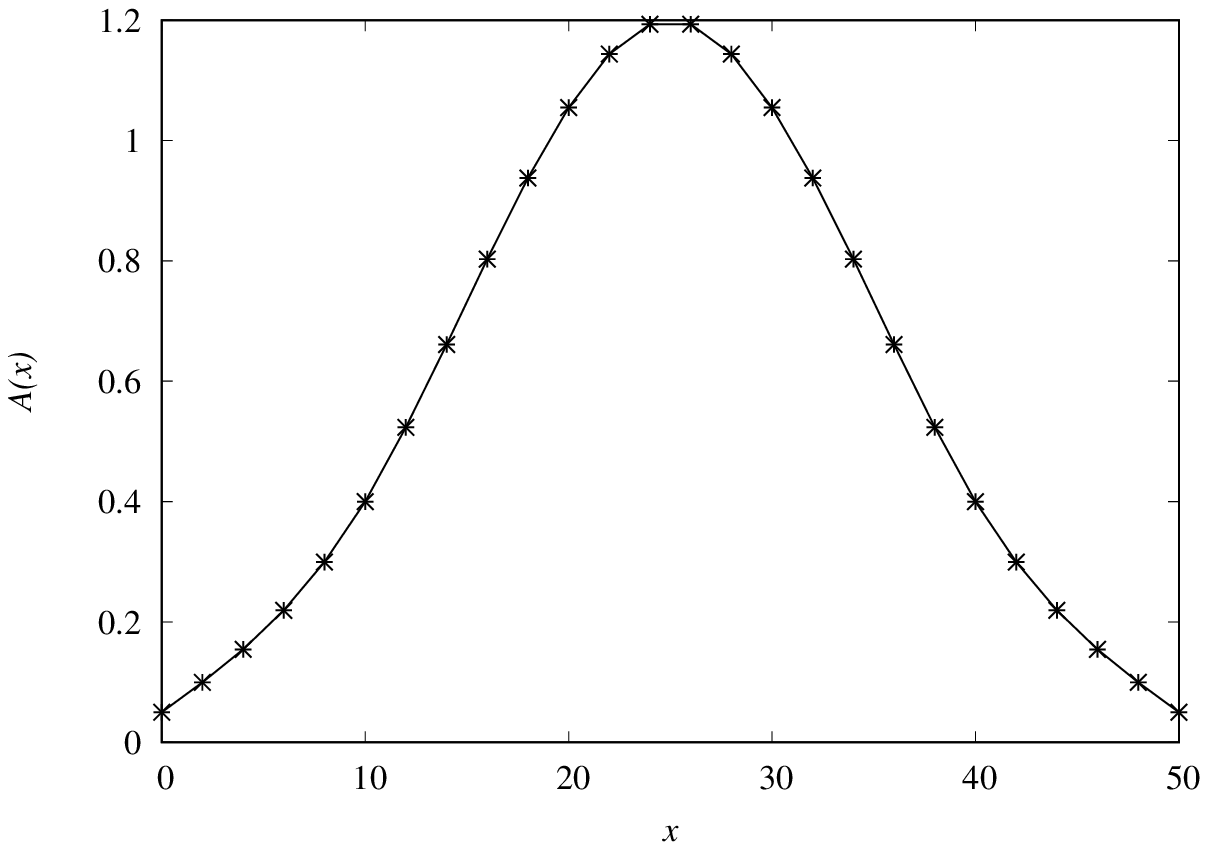}}
\end{minipage}
\hfill
\begin{minipage}[h!]{0.52\linewidth}
\center{\includegraphics[width=\linewidth]{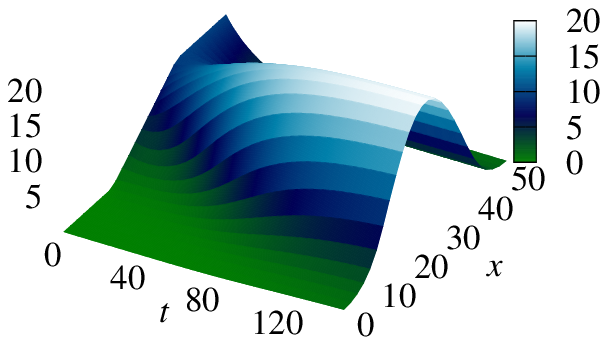}}
\end{minipage}
\vspace*{-5mm}
\caption{The space-dependent technological term $A(x)$ (left) and the solution of the direct problem~(\ref{eq:spartial_Solow_m}) for space-dependent technological level $A(x)$ and $q_{ex}(k)$ (right).}
\label{ris:inidata_A(x)}
\end{figure}

We consider a space-dependent technological level $A(x)$ demonstrated at figure~\ref{ris:inidata_A(x)} (left). Then the solution of the direct problem~(\ref{eq:spartial_Solow_m}) for the exact function $q_{ex}$ demonstrates on figure~\ref{ris:inidata_A(x)} (right).

\begin{figure}[t]
\begin{center}
\includegraphics[width=1\linewidth]{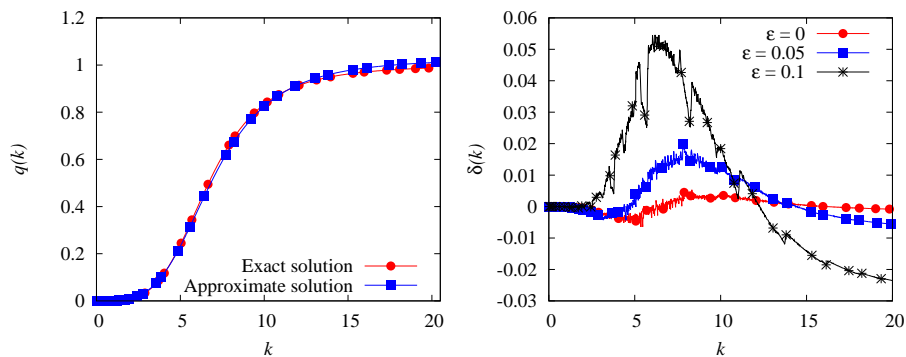}
\caption{The exact $q_{ex}(k)$ and approximate $q_\varepsilon(k)$ solutions of inverse problem for spatial Solow model for error level in data~(\ref{eq:ip_data}) $\varepsilon=0.1$ (left) and difference  $\delta (k)$ of exact and approximate solutions of inverse problem for different noise levels $\varepsilon=0, 0.05, 0.1$ in measurements (right) for $M=5$, $N=6$ for space-dependent technological level $A(x)$.}
\label{ris:func_q_A(x)}
\end{center}
\vspace*{-2mm}
\end{figure}

Using the same simulated dataset (see Section~\ref{simulated_dataset}) the inverse problem~(\ref{eq:spartial_Solow_m}),~(\ref{eq:ip_data}) is solved for number of measurements $M=5$, $N=6$ and differents error level $\varepsilon = 0, 0.05, 0.1$. The results are collected to table~\ref{tab:tab3} and demonstrated at figure~\ref{ris:func_q_A(x)}. Note, that the results of inverse problem solution are the same as for constant technological term $A$ (see Section~\ref{constantA}), i.e. accuracy in relative error $\rho_\varepsilon$ is less than $10^{-2}$, maximum of absolute difference of exact and approximated solutions of inverse problem $\max|\delta(k)|$ is the same order of $10^{-2}$. The difference of exact and approximated solutions of inverse problem $\delta (k)$ for $\varepsilon = 0, 0.05, 0.1$ in inverse problem data~(\ref{eq:ip_data}) is plotted at figure~\ref{ris:func_q_A(x)} (right). We can see that such error in reconstruction of parameters $\alpha_1$, $\alpha_2$ and $p$ (see table~\ref{tab:tab3}) is not critical to the behavior of the function $q(k)$ (see figure~\ref{ris:func_q_A(x)} from the left that demonstrated the exact and reconstructed solutions of inverse problem for error level in data~(\ref{eq:ip_data}) $\varepsilon=0.1$ and $M=5$, $N=6$).
\begin{table}[ht!]
\caption{Reconstructed parameters in function $q_\varepsilon(k)$ for different error levels $\varepsilon = 0, 0.05, 0.1$ in measurements~(\ref{eq:ip_data}) for $M=5$, $N=6$ and space-dependent technological term $A(x)$.} \label{tab:tab3}
{\footnotesize
\begin{tabular}{ccccc} \hline
\multirow{2}{*}{Parameters} & \multirow{2}{*}{Exact values} & \multicolumn{3}{c}{Error level in inverse problem data} \\
& & $\varepsilon = 0$ & $\varepsilon = 0.05$ & $\varepsilon = 0.1$ \\
\hline
$\alpha_1$ & 0.0005 & 0.00055  & 0.00065 & 0.00018 \\
$\alpha_2$ & 0.0005 & 0.00055 & 0.00064 & 0.00018 \\
$p$ & 4 & 3.9409 & 3.843 & 4.5525 \\
\hline
$\max|\delta(k)|$ & & 0.005 & 0.02 & 0.05 \\
$\rho_\varepsilon$ & & 0.001 & 0.005 & 0.009 \\
$J(q_\varepsilon)$ & & $9\cdot 10^{-5}$ & $3.8\cdot 10^{-2}$ & 0.118 \\
\hline
\end{tabular}}
\vspace*{-2.5mm}
\end{table}

\begin{figure}[t]
\begin{center}
\includegraphics[width=1\linewidth]{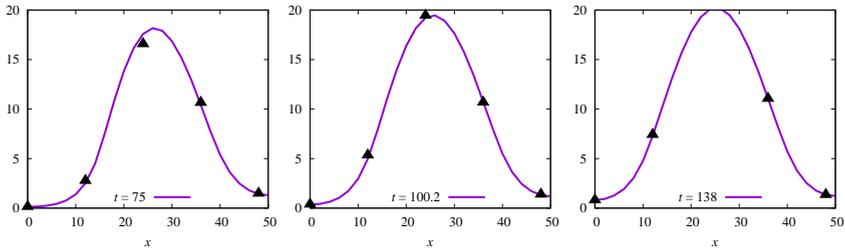}
\caption{The solution $k(x,t_j;q_\varepsilon)$ of the direct problem~(\ref{eq:spartial_Solow_m}) for reconstructed $q_\varepsilon(k)$ with error level in measurements $\varepsilon = 0.1$ and points of measurement are $M=5$, $N=6$ for space-dependent technological term. Here $t_1=75$, $t_3=100.2$, $t_6=138$.}
\label{ris:direct_pr_solution_A(x)}
\end{center}
\vspace*{-2mm}
\end{figure}

The solution $k(x,t;q_\varepsilon)$ of spatial Solow mathematical model for reconstructed $q_\varepsilon(k)$ and measured data~(\ref{eq:ip_data}) for $\varepsilon = 0.1$ is demonstrated on figure~\ref{ris:direct_pr_solution_A(x)} in case of space-dependent technological level $A(x)$. Note that capital stocks $k(x,t_j;q_\varepsilon)$, $j=1,2,3$ (purple lines) are close to the measurement points $f^\varepsilon$ (black triangles) as expected.

\section{Conclusion and outlook}
Today, economists use Solow's sources-of-growth accounting to estimate the separate effects on economic growth of technological change, capital, and labor. One important use of the Solow growth model is to estimate the share of observed growth that has resulted from growth in Total Factor Productivity (TFP), rather than from the application of increased inputs - labor, capital, and human capital (increased productive skills resulting from education and training.) Using the Solow model to approximate the output that would result in the absence of any change in TFP, you can then subtract this value from the output actually produced, and attribute the difference to TFP growth. The growth Solow model is the starting point of all analyses in modern economic growth theories, thus understanding of the model is essential to understanding the theories of the Solow growth.

The differential evolution algorithm is applied to the optimization problem of the production function $q(k)$ reconstruction for the spatial Solow model using additional measurements of GDP type for fixed space and time. Despite the fact that the considered inverse problem is ill-posed, numerical calculations show a good result with the accuracy of recovery of the production function is more than 95\% (in case of error level in measured data $10\%$). We compare the results with calculations from paper~\cite{Engbers_2014} where the authors applied Tikhonov regularization and gradient method for solving the regularized optimization problem. In the case of full measured data (that means $f^\varepsilon(x,t) = k(x,t) + \varepsilon (x,t)$) and error level $10\%$ the accuracy of reconstruction of production function was 80\% for both cases of technological levels. The reason consists in sensitivity of local regularization methods to an initial approximation while the reconstruction results for DE approach do not depend on initial population.

\bibliography{Krivorotko_reference}

\end{document}